# A strategy on prion AGAAAAGA amyloid fibril molecular modelling

Jiapu Zhang[*1]. David D.W. Liu[2]

[1] *School of Sciences, Information Technology and Engineering, University of Ballarat, Mount Helen Campus, Chancellor Drive, Mount Helen, Ballarat, VIC 3350, Australia, Emails: j.zhang@ballarat.edu.au, jiapu zhang@hotmail.com, Phone: (61)423 487 360.*

[2] *School of Medicine and Dentistry, James Cook University, Douglas Campus, University Road, Building 047, Douglas, Townsville, QLD 4814, Australia, Email: david.liu@jcu.edu.au, Phone: (61) 401 955 199.*

[*] Corresponding author

**Abstract:** X-ray crystallography and nuclear magnetic resonance (NMR) spectroscopy are two powerful tools to determine the protein 3D structure. However, not all proteins can be successfully crystallized, particularly for membrane proteins. Although NMR spectroscopy is indeed very powerful in determining the 3D structures of membrane proteins, same as X-ray crystallography, it is still very time-consuming and expensive. Under many circumstances, due to the noncrystalline and insoluble nature of some proteins, X-ray and NMR cannot be used at all. Computational approaches, however, allow us to obtain a description of the protein 3D structure at a submicroscopic level.

To the best of the authors' knowledge, there is little structural data available to date on the AGAAAAGA palindrome in the hydrophobic region (113–120) of prion proteins, which falls just within the N-terminal unstructured region (1–123) of prion proteins. Many experimental studies have shown that the AGAAAAGA region has amyloid fibril forming properties and plays an important role in prion diseases. Due to the noncrystalline and insoluble nature of the amyloid fibril, little structural data on the AGAAAAGA is available. This paper introduces a simple molecular modelling strategy to address the 3D atomic-resolution structure of prion AGAAAAGA amyloid fibrils. Atomic-resolution structures of prion AGAAAAGA amyloid fibrils got in this paper are useful for the drive to find treatments for prion diseases in the field of medicinal chemistry.

# 1 Introduction
Prion diseases are invariably fatal and highly infectious neurodegenerative diseases affecting humans and animals. The neurodegenerative diseases such as Creutzfeldt-Jakob disease (CJD), variant Creutzfeldt-Jakob diseases (vCJD), iatrogenic CreutzfeldtJakob disease (iCJD), familial Creutzfeldt-Jakob disease (fCJD), sporadic Creutzfeldt-Jakob disease (sCJD), Gerstmann-Straussler-Scheinker syndrome (GSS), Fatal Familial Insomnia (FFI), Kuru in humans, scrapie in sheep, bovine spongiform encephalopathy (BSE or mad-cow disease) in cattle, chronic wasting disease (CWD) in white-tailed deer, elk, mule deer, moose, transmissible mink encephalopathy (TME) in mink, feline spongiform encephalopathy (FSE) in cat, exotic ungulate encephalopathy (EUE) in nyala, oryx, greater kudu, and spongiform encephalopathy (SE) in ostrich etc belong to prion diseases. By now there have not been some effective therapeutic approaches or medications to treat all these prion diseases.

Prion diseases are amyloid fibril diseases. The normal cellular prion protein ($PrP^C$) is rich in α-helices but the infectious prions ($PrP^{Sc}$) are rich in β-sheets amyloid fibrils. The conversion of $PrP^C$ to $PrP^{Sc}$ is believed to involve a conformational change from a predominantly α-helical protein (about 42% α-helix and 3% β-sheet) to a protein rich in β-sheets (about 30% α-helix and 43% β-sheet) [1].

Many experimental studies such as [2, 3, 4, 5, 6, 7, 8, 9, 10, 11, 12] have shown two points: (1) the hydrophobic region (113-120) AGAAAAGA of prion proteins is critical in the conversion from a soluble $PrP^C$ into an insoluble $PrP^{Sc}$ fibrillar form; and (2) normal AGAAAAGA is an inhibitor/blocker of prion diseases. PrP lacking the palindrome could not convert to $PrP_{Sc}$ and also did not generate proteinase-K resistance. The presence of residues 119 and 120 (the two last residues within the motif AGAAAAGA) seems to be crucial for this inhibitory effect. The replacement of glycine at residues 114 and 119 by alanine led to the inability of the peptide to build fibrils but it nevertheless increased. The A117V variant is linked to the GSS disease. Etc.. Furthermore, we computationally clarified by ourselves that prion AGAAAAGA segment indeed has an amyloid fibril forming property [13, 14, 15]: Figure 1.

< Figure 1 >

*Figure 1: Prion AGAAAAGA (113-120) is surely and clearly identified as the amyloid fibril formation region, because its energy is less than the amyloid fibril formation threshold energy of -26 kcal/mol [33].*

By now the 3D structure of the AGAAAGA palindrome peptide has not been known. The physiological conditions such as pH [5] and temperature [16] will affect the propensity to form fibrils in this region. However, laboratory experiences have shown that using traditional experimental methods is very difficult to obtain atomic-resolution structures of AGAAAAGA due to the noncrystalline and insoluble nature of the amyloid fibril [17, 18] and its unstable nature. By introducing novel mathematical canonical dual formulations and computational approaches, in this paper we may construct atomic-resolution molecular structures for prion (113-120) AGAAAAGA amyloid fibrils.

The 3D atomic resolution structure of PrP (106-126), i.e. TNVKHVAGAAAAGAVVGGLGG, can be looked as the structure of a control peptide [19, 20]. Ma and Nussinov (2002) established homology structure of AGAAAAGA and its molecular dynamics simulation studies [21]. Recently, Wagoner et al. studied the structure of GAVAAAAVAG of mouse prion protein [22, 16]. All these structures (including the structures of this current paper) of amyloid cross-β spines can be reduced to the ones of 8 classes of steric zippers reported in [23]. Thus, we might be able to say that the molecular modeling strategy reported in this study should be reliable.

Many studies have indicated that computational approaches or introducing novel mathematical formulations and physical concepts into molecular biology can significantly stimulate the development of biological and medical science. Various computer computational approaches were used to address the problems related to "amyloid fibril" [24, 25, 26, 27, 28, 29]. Here, we would like to use the simulated annealing evolutionary computations to build the optimal atomic-resolution amyloid fibril models in hopes to be used for controlling prion diseases. The templates from the Protein Data Bank (pdb) (http://www.rcsb.org) and the modeling computational algorithms are different among this current study and previous studies [14, 15]. Zhang (2011) used all the pdb files of [23] and the Insight II package to establish the models [14]. Zhang et al. (2011) used 3FVA.pdb (the NNQNTF (173-178) segment from elk prion protein) and hybrid simulated annealing discrete gradient (SADG) method to build the models [15]. This current study uses 3NHC.pdb (the GYMLGS (127-132) segment from human M129 prion protein) as the modelling template and simulated annealing evolutionary computations (SAEC) to build the models, where SAECs were got from the hybrid algorithms of [30] by simply replacing the DG method by the SA algorithm of [31, 15] and numerical computational results show that SAECs can successfully pass the test of more than 40 well-known benchmark global optimization problems.

The atomic structures of all amyloid fibrils revealed steric zippers, with strong van der Waals (vdw) interactions between β-sheets and hydrogen bonds (HBs) to maintain the β-strands [23]. The vdw contacts of atoms are described by the Lennard-Jones (LJ) potential energy:

$$V_{LJ}(r) = 4\varepsilon [(\sigma/r)^{12} - (\sigma/r)^{6}], \quad (Eq.\ 1)$$

where ε is the depth of the potential well and σ is the atom diameter; these parameters can be fitted to reproduce experimental data or deduced from results of accurate quantum chemistry calculations. The $(\sigma/r)^{12}$ term describes repulsion and the $-(\sigma/r)^{6}$ term describes attraction. If we introduce the coordinates of the atoms whose number is denoted by N and let ε=σ=1 be the reduced units, the formula (Eq. 1) becomes

$$f(x) = 4 \sum_{i=1}^{N} \sum_{j=1, j<i}^{N} (1/t_{ij}^{6} - 1/t_{ij}^{3}), \quad (Eq.\ 2)$$

where $t_{ij} = (x_{3i-2} - x_{3j-2})^2 + (x_{3i-1} - x_{3j-1})^2 + (x_{3i} - x_{3j})^2$, $(x_{3i-2}, x_{3i-1}, x_{3i})$ is the coordinates of atom i, N≥2. The minimization of LJ potential f(x) on $R^n$ (where n = 3N) is an optimization problem:

$$\min f(x) \text{ subject to } x \in R^{3N}. \quad (Eq.\ 3)$$

Similarly as (Eq. 1) – the potential energy for the vdw interactions between β-sheets:

$$V_{LJ}(r) = A/r^{12} - B/r^6, \quad (Eq.\ 4)$$

the potential energy for the HBs between the β-strands has the formula

$$V_{HB}(r) = C/r^{12} - D/r^{10}, \quad (Eq.\ 5)$$

where A,B,C,D are given constants. Thus, the amyloid fibril molecular modelling problem is deduced into well solve the mathematical optimization problem (Eq. 3).

This paper is organized as follows. In Section 2, we first describe how to build the prion AGAAAAGA amyloid fibril molecular models, and then explain how the models with only 6 variables are built and can be solved by any optimization algorithm. At the end of Section 2 the models are done a little refinement by Amber 11 [32]. At last, we conclude that when using the time-consuming and costly X-ray crystallography or NMR spectroscopy we still cannot determine the protein 3D structure, we may introduce computational approaches or novel mathematical formulations and physical concepts into molecular biology to study molecular structures. This concluding remark will be made in the last section.

## 2 Prion AGAAAAGA amyloid fibril models' Molecular Modeling and Optimizing

Constructions of the AGAAAAGA amyloid fibril molecular structures of prion 113–120 region are based on the most recently released experimental molecular structures of human M129 prion peptide 127–132 (PDB entry 3NHC released into Protein Data Bank on 04-AUG-2010). The atomic-resolution structure of this peptide is a steric zipper, with strong vdw interactions between β-sheets and HBs to maintain the β-strands (Figure 2).

< Figure 2 >

*Figure 2: Protein fibril structure of human M129 prion GYMLGS (127–132) (PDB ID: 3NHC). The dashed lines denote the hydrogen bonds. A, B, ..., K, L denote the chains of the fibril.*

In Figure 2 we see that G (H) chains (i.e. β-sheet 2) of 3NHC.pdb can be obtained from A (B) chains (i.e. β-sheet 1) by

$$G(H) = ((1\ 0\ 0)^T, (0\ -1\ 0)^T, (0\ 0\ -1)^T)\ A(B) + (9.07500\ 4.77650\ 0.00000)^T \quad (Eq.\ 6)$$

and other chains can be got by

$$I(J) = G(H) + (0\ 9.5530\ 0)^T,\ K(L) = G(H) + (0\ -9.5530\ 0)^T, \quad (Eq.\ 7)$$

$$C(D) = A(B) + (0\ 9.5530\ 0)^T,\ E(F) = A(B) + (0\ -9.5530\ 0)^T. \quad (Eq.\ 8)$$

Basing on the template 3NHC.pdb from the Protein Data Bank, three prion AGAAAAGA palindrome amyloid fibril models - an AAAAGA model (Model 1), a GAAAAG model (Model 2), and an AAAAGA model (Model 3) - will be successfully constructed in this paper. AB chains of Models 1-3 were respectively got from AB chains of 3NHC.pdb using the mutate module of the free package Swiss-PdbViewer (SPDBV Version 4.01) (http://spdbv.vital-it.ch). It is pleasant to see that almost all the hydrogen bonds are still kept after the mutations; thus we just need to consider the vdw contacts only. Making mutations for GH chains of 3NHC.pdb, we can get the GH chains of Models 1-3. However, the vdw contacts between A chain and G chain, between B chain and H chain are too far at this moment (Figures 3-5).

< Figure 3 >
*Figure 3: At initial state, the vdw contacts between AG chains (_-sheet 1) and BH chains (β-sheet 2) of Model 1 are very far.*

< Figure 4 >
*Figure 4: At initial state, the vdw contacts between AG chains (_-sheet 1) and BH chains (β-sheet 2) of Model 2 are very far.*

< Figure 5 >
*Figure 5: At initial state, the vdw contacts between AG chains (_-sheet 1) and BH chains (β-sheet 2) of Model 3 are very far.*

Seeing Figures 3-5, we may know that for Models 1-3 at least two vdw interactions between A.ALA3.CB-G.ALA4.CB, B.ALA4.CB-H.ALA3.CB should be maintained. Fixing the coordinates of A.ALA3.CB and B.ALA4.CB, letting the coordinates of G.ALA4.CB and H.ALA3.CB be variables, we may get a simple LJ potential energy minimization problem (Eq. 3) just with six variables. For solving this six variable optimization problem, any optimization computational algorithm can be used to solve this low-dimensional mathematical optimization problem; for example, in this paper we may use the SAEC algorithm. Setting the coordinates of G.ALA4.CB and H.ALA3.CB as initial solutions, running a SAEC algorithm, for Models 1-3 we get

$$G(H) = ((1\ 0\ 0)^T,\ (0\ -1\ 0)^T,\ (0\ 0\ -1)^T)\ A(B) + (-0.703968\ 7.43502\ -0.33248)^T.\quad (Eq.\ 9)$$

By (Eq. 9) we can get close vdw contacts between A chain and G chain, between B chain and H chain (Figures 6-8).

< Figure 6 >
*Figure 6: After LJ potential energy minimization, the vdw contacts between AG chains (β-sheet 1) and BH chains (β-sheet 2) of Model 1 become very closer.*

< Figure 7 >
*Figure 7: After LJ potential energy minimization, the vdw contacts between AG chains (β-sheet 1) and BH chains (β-sheet 2) of Model 2 become very closer.*

< Figure 8 >
*Figure 8: After LJ potential energy minimization, the vdw contacts between AG chains (β-sheet 1) and BH chains (β-sheet 2) of Model 3 become very closer.*

Furthermore, we may employ the Amber 11 package [32] to slightly optimize Models 1-3 and at last get Models 1-3 with stable total potential energies (Figures 9-11). The other CDIJ and EFKL chains can be got by parallelizing ABGH chains in the use of mathematical formulas (Eq. 7)-(Eq. 8).

< Figure 9 >

*Figure 9: After Amber 11 refinement of the whole potential energy minimization, the optimal structure of prion AGAAAA amyloid fibril Models 1.*

< Figure 10 >

*Figure 10: After Amber 11 refinement of the whole potential energy minimization, the optimal structure of prion GAAAAG amyloid fibril Models 2.*

< Figure 11 >

*Figure 11: After Amber 11 refinement of the whole potential energy minimization, the optimal structure of prion AAAAGA amyloid fibril Models 3.*

# 3 Conclusion

X-ray crystallography is a powerful tool to determine the protein 3D structure. However, it is time-consuming and expensive, and not all proteins can be successfully crystallized, particularly for membrane proteins. Although NMR spectroscopy is indeed a very powerful tool in determining the 3D structures of membrane proteins, it is also time-consuming and costly. Due to the noncrystalline and insoluble nature of the amyloid fibril, little structural data on the prion AGAAAAGA segment is available. Under these circumstances, the novel simple strategy introduced in this paper can well do the molecular modeling of prion AGAAAAGA amyloid fibrils. This indicated that computational approaches or introducing novel mathematical formulations and physical concepts into molecular biology can significantly stimulate the development of biological and medical science. The optimal atomic-resolution structures of prion AGAAAAGA amyloid fibils presented in this paper are useful for the drive to find treatments for prion diseases in the field of medicinal chemistry.

# Acknowledgments

This research was supported by a Victorian Life Sciences Computation Initiative (VLSCI) grant number VR0063 on its Peak Computing Facility at the University of Melbourne, an initiative of the Victorian Government. The authors appreciate the anonymous referees for the numerous insightful comments, which have improved this paper greatly.

Figures 1-11 can be found in the pdf file**Cent Eur J Biol - A strategy on prion AGAAAAGA amyloid fibril molecular modeling - 3NHC.pdf** in website http://blog.51xuewen.com/jiapu_zhang_PhD_MSc_MSc_BSc/article_42219.htm